\documentclass[12pt]{article} 
\usepackage[sectionbib]{natbib}
\usepackage{array,epsfig,fancyheadings,rotating,xcolor}
\usepackage[]{hyperref}  
\usepackage{sectsty, secdot,todonotes}
\sectionfont{\fontsize{12}{14pt plus.8pt minus .6pt}\selectfont}
\subsectionfont{\fontsize{12}{14pt plus.8pt minus .6pt}\selectfont}

\usepackage{amsmath}
\usepackage{amssymb}
\usepackage{amsfonts}
\usepackage{multirow}
\usepackage{amsthm}
\usepackage{soul}

\newtheorem{theorem}{Theorem}
\newtheorem{lemma}{Lemma}
\newtheorem{corollary}{Corollary}
\newtheorem{proposition}{Proposition}
\theoremstyle{definition}
\newtheorem{definition}{Definition}

\newtheorem{remark}{Remark}
\newcommand{\eps}{\varepsilon}

\newcommand{\EXP}{\mathbb{E}}
\newcommand{\prob}{\mathbb{P}}
\newcommand{\normal}{\mathbf{n}}
\newcommand{\gen}{\mathcal{L}}
\newcommand{\R}{\mathbb{R}}

\begin{document}
	

\begin{center}
	\Large {\textbf	Level sets and drift estimation for \\ reflected Brownian motion with drift}
\end{center}

\begin{center}
	Alejandro Cholaquidis$^*$, Ricardo Fraiman$^*$,\\
	Ernesto Mordecki$^*$, Cecilia Papalardo$^*$\\
	$^*$ Universidad de la Rep\'ublica
\end{center}

\begin{abstract}
We consider the estimation of  the drift and the level sets of the stationary distribution 
of a Brownian motion with drift, reflected in the boundary of a compact set $S\subset\R^d$, 
departing from the observation of a  trajectory  of this process. 
We obtain the uniform consistency and rates of convergence for the proposed kernel based estimators. 
This problem has relevant applications in ecology, 
in estimating the home-range and the core-area of an animal based on tracking data.  
Recently, the problem of estimating the domain of a reflected Brownian motion
was considered in \cite{ch:16},
in this case the stationary distribution is uniform and the estimation of the core-area, 
defined as a level set of the stationary distribution, is meaningless. 
We also give an estimator of the drift function, based on the increments of the process.
In order to prove our results, 
some new theoretical properties of the reflected Brownian motion with drift are obtained,
under fairly general assumptions. 
These properties allow us to perform the estimation for flexible regions close to reality. 
The theoretical findings are illustrated in simulated and real data examples.
\end{abstract}

\section{Introduction}\label{sec:intro}

Given a Reflected Brownian
Motion with Drift (RMBD) inside a (smooth enough) compact domain S, we will consider three statistical problems.
The first one is the estimation of the the density of the stationary distribution. The second one is the
estimation of the level sets of this density, with and without shape restrictions.
Lastly, we consider the problem of the estimation of the drift function. 
The practical motivation of such problems will be made explicit in
the following paragraphs.

Level set estimation can be placed into the field of non-parametric set estimation, where the goal
is reconstructing (in the statistical sense) an unknown set $S$ from random data related
to $S$. Usually, such random information comes from a sample of independent points
drawn form an absolutely continuous distribution with density $f$ , and the target set
is either the support of $f$ or a level set of 
the type $\{x : f (x) > \lambda\}$ (which, depending on $\lambda$,
can be seen as a sort of “substantial support” of the underlying distribution).

There are, however, two important practical applications of set estimation techniques in
which the assumption of independence is clearly unsuitable and the above mentioned
RBMD approach might be particularly well-motivated. These are the problems of
estimation of the so-called home range (the region where an individual of an animal
species develops its activities, Burt, 1943) and that of the core area (the sub-region of
the home range where the individual spends most of its time; see \cite{hayne:49}, \cite{worton:87}). 
Recent advances in animal tracking technology allows an almost continuous record of the movement.
Therefore in both cases it might be reasonable to assume that the sample information
comes from a grid of points taken along the (random) trajectory followed by the animal
during its activities. In our setup, core-areas can be modelled by the level sets of the stationary distribution, while the drift function provides information about the dynamics of the movement of the animal.

Such point of view has been followed in \cite{ch:16} where the home range
$S$ is identified with the support of the stationary distribution of a Reflected
Brownian Motion (with no drift) and a suitable set estimator is proposed and analyzed.
However, under the (quite natural) regularity conditions on S considered in that paper
the stationary distribution of the reflected Brownian motion is necessarily uniform, (\cite{burdzy:06}).
This is somewhat restrictive under different points of view; in particular, the problem
of core area estimation (which should be addressed in terms of level set estimation) is
meaningless or uninteresting.

Thus, the main contribution of the present paper is to extend the approach in 
\cite{ch:16}
to the case of Reflected Brownian Motion with Drift. Such extension is indeed
substantial and far beyond a simple technical generalization, as it allow us to address
the estimation of the core area in terms of the estimation of an appropriate
level set of the stationary distribution of a RBMD in the home range $S$. Note
that this can be done since such stationary distribution is non-uniform, in general.

As a by-product, the paper provides explicit conditions for the existence and geometric ergodicity 
of a RBMD on the domain, The drift estimation problem is also addresses. 

With a different approach, an exponential rate 
in the estimation of the stationary distribution
has been also obtained for ergodic diffusions in unbounded domains (see \cite{dal:05}). See also \cite{clr:16}, where a similar problem is considered. 
The estimation of the stationary distribution of a stochastic differential equation with drift, 
without reflection, 
has been studied by several authors, 
some earlier results are given in \cite{ver99} who adapts the results of \cite{cl86} to this setting. 
More recently, \cite{dal:04} estimates the drift and the stationary distribution for the same model but without reflection,  
whereas \cite{go:04} considers estimation problems for one dimensional diffusions with and without reflection.

Before introducing the formal framework, we
discuss briefly the application of the proposed method to the estimation from animal tracking data of the  core-area and drift. 
For a description of home-range estimation see for instance \cite{ch:16} and the references therein.

\subsection{Roadmap}
This paper is organized as follows.
In Section \ref{setup}, we discuss conditions for the existence, uniqueness and geometric ergodicity of the reflected Brownian motion 
with drift. The main results in this section are given in Propositions \ref{prop0} and  \ref{properg}. 
Proposition \ref{prop0} gives sufficient conditions for Harris recurrence and for the the domain to be non-trap for the RBMD process $\{X_t\}$ 
(a condition introduced in \cite{burdzy:06}, which we describe in Section  \ref{setup}).  

In Proposition  \ref{properg} we show that if the domain is non-trap, 
we have an exponential rate of convergence to the stationary distribution for the total variation norm.  All the proofs for this section are given in Appendix B.
In Section \ref{driftest} we obtain in Theorem \ref{teoest1} strong uniform convergence rates for kernel estimators of the stationary distribution based on a trajectory of the RBMD. In Corollary \ref{corlevset} and Theorem \ref{teoest2} we prove the strong consistency of two different families of level sets estimators with respect to the Hausdorff distance. 
	The case of the estimation of level sets with a given content is considered in Theorem \ref{thgtau}. 
We also derive in Theorem \ref{condrift} consistent estimators of the drift function. 
Lastly, in Section \ref{examples} we consider some simulated and real data examples to illustrate the behaviour of the estimation methods described in the paper.

\section{Reflected Brownian motion with drift}\label{setup}

In this section we establish conditions for the existence of a reflected Brownian motion with drift and its stationary distribution, and study the connections between these conditions and some geometric constraints on its support. 
\subsection{Notation}
Given a set $S\subset\mathbb{R}^d$, we denote by $\partial S$, $\textnormal{int}(S)$, and $\overline{S}$ the boundary, interior, and closure of $S$, respectively. If $S$ is a finite set, we denote by $\# S$ its cardinal.
The Borel sigma algebra in $S$ will be denoted by $\mathcal{B}(S)$. 
We denote by $\left\langle\cdot,\cdot \right\rangle$  the usual inner product in $\mathbb{R}^d$ 
and by $\left\|\cdot\right\|$ the Euclidean norm.
The closed ball of radius $\eps$ centred at $x$ is denoted by $\mathcal{B}(x,\eps)$, while the open ball is denoted by $\mathring{\mathcal{B}}(x,\eps)$.  
Given $\eps>0$ and a bounded set $A\subset \mathbb{R}^d$, $B(A,\epsilon)$ denotes the parallel set $B(A,\epsilon)=\{x\in \mathbb{R}^d\colon\ d(x,A)\leq \epsilon\}$ where $d(x,A)=\inf\{\|x-a\|\colon\ a\in A\}$. The $d$-dimensional Lebesgue measure on $\mathbb{R}^d$ will be denoted by $\mu_L$.
\subsection{An implicit definition of the Reflected Brownian Motion}
Let $D$ be a bounded domain in $\mathbb{R}^d$ (that is, a bounded connected open set) such that $\partial D$ is $C^2$. 
Given a $d$-dimensional Brownian motion $\{B_t\}_{t\geq 0}$ departing from 
$B_0=0$ 
defined on a filtered probability space $(\Omega,\mathcal{F},\{\mathcal{F}_t\}_{t\geq 0},\prob_x)$, 
we are concerned with the problem of the existence and uniqueness of the solution of a reflected stochastic differential equation on $\overline{D}$ given by
\begin{equation}\label{sde}
X_t= X_0+ B_t+ \int_0^t\mu(X_s)ds+\int_0^t\normal(X_s)\xi(ds),
\quad\text{ where } X_t\in \overline{D},\ \forall t\geq 0.
\end{equation} 
whose drift, $\mu(x)$, is assumed to be Lipschitz, while $\normal(x)$ denotes the inner unit vector at the boundary point $x\in\partial D$;  
this boundary satisfies some regularity conditions (to be specified later).
This equation is called a \emph{Skorokhod stochastic differential equation}.
Its solution is a pair of stochastic processes $\{X_t,\xi_t\}_{t\geq 0}$, 
the first coordinate $\{X_t\}_{t\geq 0}$ is a \emph{reflected diffusion}, which we call a 
\emph{reflected Brownian motion with drift} (RBMD), 
and $\{\xi_t\}_{t\geq 0}$ is the corresponding \emph{local time}, that is, a one-dimensional continuous non-decreasing process with $\xi_0=0$ that satisfies
\[
\xi_t=\int_0^t\mathbb{I}_{\{X_s\in\partial D\}}d\xi_s.
\]
Since we have assumed that $\partial D$ is $C^2$, we know that a ball of positive radius rolls freely inside and outside $\overline{D}$ (see \cite{walther:99}). Then, 
by using the same arguments used to prove proposition 3 in \cite{ch:16}, 
we can ensure that the geometric 
shape
conditions for the existence and uniqueness of a solution of equation \eqref{sde}, as required in \cite{saisho:87}, are satisfied. 
From Theorem 5.1 in \cite{saisho:87} it follows that there exists a unique strong solution of the Skorokhod stochastic differential equation \eqref{sde}. The solution is a strong solution in the sense of definition 1.6 in \cite{ikeda}.

\begin{remark} \label{rem}
   There exists a unique positive function $p(s,x,t,y)$ satisfying	$\mathbb{P}(X_t\in \Gamma|X_s=x)=P(s,x,t,\Gamma)=\int_{\Gamma}p(s,x,t,y)dy$ and, by theorem 3.2.1 of \cite{sv:97}, 
	the function $p$ satisfies the forward equation
	$\partial_s p+\gen^* p=0$ and $\lim_{s\rightarrow t^-} p(s,.,t,y)=\delta_y$, where $\delta_y$ is the point-mass at $y$ and $\mathcal{L}^*$ is the adjoint of $\mathcal{L}$, that is, $\gen^* h=\frac{1}{2}\Delta h-\langle \mu , \nabla h \rangle.$
\end{remark}

\subsection{Ergodic properties}

We now introduce the notions of invariant measure and ergodic process, following \cite{MTb}. 
\begin{definition} 
A probability measure $\pi$ on $S$ is said to be an \emph{invariant measure} for a time-homogeneous Markov process $\{Z_t\}_{t\geq 0}$ if 
$\int_S \prob_x(Z_t\in A)\pi(dx)=\pi(A)$, for all $t>0$ and all $A\in \mathcal{B}(S)$. \end{definition}


\begin{definition}  A Markov process $\{Z_t\}_{t\geq 0}$ with state space $S$ is \emph{ergodic} if there exists an invariant probability measure $\pi$ such that 
$\lim_{t\rightarrow +\infty} \big\|\prob_x(Z_t\in \cdot)-\pi(\cdot)\big\|_{TV}=0\quad \forall x\in S$. Here $\Vert  \mu \Vert_{TV}$ stands for the total variation norm of the measure $\mu$. 
In this case $\pi$ is called a \emph{stationary distribution}.	
\end{definition}

\begin{remark} \label{green}
	If the drift is given by the gradient of some function $f$, i.e, $\mu(x)=\frac{1}{2}\nabla f(x)$, by Green's formula,  there exists a unique stationary distribution and is given by $\pi(dx)=ce^{-f(x)}\mathbb{I}_Ddx=g(x)dx$, where $c$ is the normalization constant. 
\end{remark}
\begin{definition} \label{gergo} A Markov process $\{Z_n\}_{n\in \mathbb{N}}$ with state space $S$ is called geometrically ergodic if there exists an invariant probability $\pi$ and real numbers $0<\rho<1$ and $\gamma>0$ such that
		\begin{equation}\label{ge}
		\big|\mathbb{P}_x(Z_n\in B)-\pi(B)\big|\leq \gamma \rho^n\quad\text{for all $x\in S$ and all $B\in \mathcal{B}(S)$}.
		\end{equation}
	\end{definition}

\subsection{Harris recurrence and the trap condition.} 

Let $D\subset\mathbb{R}^d$ be an open bounded set and $\mathcal{B}\subset D$. 
Consider the first hitting time of $\mathcal{B}$ by a stochastic process $\{Z_t\}_{t\geq 0}$ 
defined by 
$T_\mathcal{B} = \inf\{t>0\colon Z_t \in \mathcal{B}\}$.  

\begin{definition} \label{hrec} A Markov process $\{Z_t\}_{t\geq 0}$ is called \emph{Harris recurrent} 
if for some $\sigma$-finite measure $\mu$, 
we have $\prob_x(T_A<\infty)=1$ whenever $\mu(A)>0$, $A\in\mathcal{B}(\overline{D})$.
\end{definition}	

Under Harris recurrence there exists a unique (up to a multiplicative constant) invariant measure 
(see \cite{akr:67}). For the RBMD we prove in Proposition \ref{prop0} a sufficient condition for Harris recurrence (taking in Definition  \ref{hrec} $\mu$ as the Lebesgue measure restricted to $D$), slightly stronger than the  \emph{non-trap} condition introduced in \cite{burdzy:06}).

\begin{definition}\label{def:trap}
	 We say that $D$ is a \emph{trap domain} 
	for the stochastic process $\{Z_t\}_{t\geq 0}$
	if 
	there exists a closed ball  $\mathcal{B}\subset D$ with positive radius
	such that $\sup_{x\in D}\mathbb{E}_xT_\mathcal{B}=\infty$, where $\mathbb{E}_x$ denotes the expectation w.r.t. $\prob_x$.  
	Otherwise $D$ is called a \emph{non-trap domain}. 
\end{definition}
The non--trap condition is mandatory to estimate the stationary distribution and the drift function, in order to visit infinitely many often a small ball at each point $x$.

It is proved in lemma 3.2 in \cite{burdzy:06} that if $\{X_t\}_{t\geq 0}$ is a reflected Brownian motion (without drift) in a connected open set $D$ with finite volume and $\mathcal{B}_1$, $\mathcal{B}_2$ are closed non-degenerate balls in $D$, then
$\sup_{x\in D}\mathbb{E}_xT_{\mathcal{B}_1}<\infty$ if and only if $\sup_{x\in D}\mathbb{E}_xT_{\mathcal{B}_2}<\infty$.

\begin{proposition} \label{prop0} Let $D \subset \mathbb R^d$ be a bounded domain such that $\partial D$ is $C^2$. Let $\{X_t\}_{t\geq 0}$ be the solution of \eqref{sde}, then for all Borel set $A$ such that $\mu_L(A\cap D)>0$,  we have that
		\begin{equation} \label{trapcond2}
	\sup_{x\in D}\mathbb{E}_x(T_{A})<\infty,
	\end{equation}
	where $\mathbb{E}_x$ denotes the expectation w.r.t. $\prob_x$, which implies Harris recurrence.  
	\end{proposition}	

The following proposition (whose proof is given in Appendix B) states that under the non-trap condition the process is geometrically ergodic.  This result can also be nicely derived using functional inequalities as has been proposed in \cite{clr:16}, see Section 3.1. 

\begin{proposition} \label{properg} 
Let $D \subset \mathbb R^d$ be a bounded domain such that $\partial D$ is $C^2$. 
Denote by $\pi$ the invariant distribution of $\{X_t\}_{t\geq 0}$.
If $D$ is a non-trap domain for $\{X_t\}_{t\geq 0}$, then 
there exist positive constants $\alpha$ and $\beta$ such that 
\begin{equation*}
\sup_{x \in D} \big\Vert \mathbb P_x (X_t \in \cdot ) -\pi(\cdot)\big\Vert_{TV} \leq \beta e^{-\alpha t}.
\end{equation*}
\end{proposition}

\section{Estimation of the drift and stationary distribution} \label{driftest}

In this section we first obtain in Theorem \ref{teoest1} strong  uniform convergence rates for the classical  kernel density estimator $\hat{g}_n$ of the density $g$ of the stationary distribution of a geometrically ergodic Markov chain. This allows to estimate the density, $g$, of the stationary distribution of the RBMD $\{X_t\}_{t\geq 0}$ by considering a sequence $\{X_{kn_1}\}_{k\in \mathbb{N}}$ (the choice of $n_1$ will be given explicitly in the proof of Proposition \ref{properg}). As it is well known, uniform convergence is crucial to obtain the convergence of level sets (see Theorem \ref{teoest2}).
	Next, in Corollary \ref{corlevset} and Theorem \ref{teoest2} we  show the convergence of two families of estimator of the level sets. We consider the case of the estimation of level sets with a given content in Theorem \ref{thgtau}. 

The proof of Theorem \ref{teoest1} is based on some ideas in \cite{cd:05},

the main difference 
being that we aim to obtain uniform convergence, to be able to estimate the level sets. In order to do so, 
we  introduce some notation.
    
	Let $\{X_n\}_{n\in \mathbb{N}}$ be a Markov process with state space $S\subset \mathbb{R}^d$ and let $\mu_0(dy)$ be an arbitrary initial distribution. Let $\mu_n(dy)$ denote the distribution of $X_n$, that is,
	$$\mathbb{P}_{\mu_0}(X_n\in A)=\int_A \mu_n(dy)\quad \forall A\in \mathcal{B}(S),$$
	where $\mathbb{P}_{\mu_0}$ indicates that the initial distribution is $\mu_0$. 
	Similarly,  $\EXP_{\mu_0}$ indicates the corresponding expectation.

	Let $K\colon \mathbb R^d \to \mathbb R$ be a bounded function such that $K\geq 0$ and $\int K(t)dt=1$. Consider the classical kernel estimator $\hat{g}_n$ based on $\{X_1,\dots,X_n\}$, 	given by
	\begin{equation*}
	\hat{g}_n(x)=\frac{1}{nh_n^d}\sum_{i=1}^nK\left(\frac{x-X_i}{h_n}\right)=\frac{1}{n}\sum_{i=1}^{n}K_h(x-y),
	\end{equation*}
	 where $h=h_n\rightarrow 0$ and $K_h(x)=K(x/h)/h^d$.\\
	The following generalization of the Bernstein inequality obtained in \cite{cd:84}, will be useful throughout the present discussion. Some sharper bounds were obtained more recently (see for instance \cite{dn07}). However the same rates of convergence are obtained from Collomb's inequality. Recall that a stochastic process $\{X_k\}_{k\in \mathbb{Z}}$ is $\varphi$-mixing if $\sup_{j\in \mathbb{Z}} \sigma(\mathcal{F}_{-\infty}^j,\mathcal{F}_{j+n}^\infty)\rightarrow 0$ as $n\rightarrow \infty$, where $\mathcal{F}_{j}^k=\sigma(X_s,j\leq s\leq k)$.

	\begin{lemma}  (Bernstein inequality for $\varphi$-mixing processes). Let $Y_i$ be a sequence of $\varphi$-mixing random variables such that $\EXP(Y_i)=0$, $|Y_i|\leq C_1$, $\EXP|Y_i|\leq \eta$, and $\EXP(Y_i^2)\leq D$. Write $\tilde{\varphi}(m)=\varphi(1)+\dots+\varphi(m)$ for each $m\in \mathbb{N}$. Then, for each $\eps>0$ and $n\in \mathbb{N}$, we have
		\begin{equation} \label{bern}
		\mathbb{P}\left(\left|\sum_{k=1}^n Y_k\right|>\eps\right)\leq 2\exp\Big(3e^{1/2}n\frac{\varphi(m)}{m}-\alpha\eps  +\alpha^2nC_2\Big),
		\end{equation}
where 
$C_2=6(D+4\eta C_1\tilde{\varphi}(m))$ 
and $\alpha$, $m$ are respectively any positive real number and any positive integer less than or equal to 
$n$ and satisfying $\alpha mC_1\leq 1/4$. The numbers $\alpha$ and $m$  may also depend on $n$. 
\end{lemma}

	\begin{theorem} \label{teoest1} Let $S\subset\mathbb{R}^d$ be a compact set and $\{X_n\}_{n\in \mathbb{N}}$ a geometrically ergodic Markov chain with state space $S$ and constants $\gamma$ and $\rho$ given by \eqref{ge}, 
	whose stationary distribution, $\pi$, has a Lipschitz density $g$ w.r.t. to Lebesgue measure. Denote by $g_1=\max_{x\in S}g(x)$ and by $C_g$ the Lipschitz constant of $g$. Let
		$\hat{g}_n(x)=\frac{1}{n}\sum_{i=1}^nK_h(x-X_i)$
		with $K:\mathbb{R}^d\rightarrow \mathbb{R}$ a non-negative bounded Lipschitz function such that $\int K(t)dt=1$ and  $\kappa=\int |u|K(u)du<\infty$. Denote by $k_1=\max K(x)$. Let $h=h_n\rightarrow 0$, $\alpha_n\rightarrow 0$, and $\beta_n\rightarrow \infty$ such that, $\beta_n h_n\rightarrow 0$, $\alpha_n=o(1/\beta_n)$, and $\log(n)/\beta_n\rightarrow 0$. Then, for all $\epsilon>0$, and for all $n>n_1$ ($n_1$ will be given in the proof), we have
		\begin{equation} \label{unifconvtasa}
		\mathbb{P}\Big(\beta_n\sup_{x\in S}\big|\hat{g}_n(x)-g(x)\big|>\epsilon\Big)\leq C \frac{\beta_n}{nh^d}+C'h\beta_n+3c\exp\Big(-\epsilon\alpha_n\frac{nh^d}{4\beta_n}-(d(d+2))\log(h)\Big),
		\end{equation}
		where $C=2k_1 \gamma \sum_{n=1}^\infty \rho^n$,  $C'=\kappa C_g$ and $c$ is a constant, depending only on $d$ and $\mu_L(S)$. 	\\
	   Moreover, if $\beta_n$ and $h_n$ fulfils also that  $\alpha_n nh^d/(\beta_n \log(n))\rightarrow \infty$, then  $\beta_n\sup_{x\in S}\big|\hat{g}_n(x)-g(x)\big|\rightarrow 0 \quad a.s.$
	   
		\end{theorem}

\begin{remark} \label{remconf}
\begin{itemize}
	\item[i)]  Taking $h=n^{-1/\nu}$ and $\beta_n=n^\gamma$, then the best attainable rate that can be derived from Theorem \ref{teoest1} is for $\gamma=\frac{1}{d+2}$, i.e, $\beta_n=\mathcal{O}(n^{1/(d+2)})$.
	\item[ii)] If we only want uniform convergence, the conditions in $h_n$ can be relaxed, and replaced by $h=\mathcal{O}((1/n)^{1/(d+1)})$. 
	\end{itemize}
\end{remark}

Using now Theorem 2 of \cite{cgmrc:06} we get the following direct corollary, which establishes the rate for  the consistency in Hausdorff distance of the boundary of the estimated level sets $\partial G_{\hat{g}_n}(\lambda)$ (where $G_g(\beta)=\{x\colon g(x)>\beta\}$). Recall that given two non-empty compact sets $A,C\subset \mathbb{R}^d$, the Hausdorff distance between $A$ and $C$ is defined as 
\begin{equation}
d_H(A,C)=\max\Big\{\max_{a\in A}d(a,C), \ \max_{c\in C}d(c,A)\Big\}, \text{ where }
	 d(a,C)=\inf_{c\in C} d(a,c).
\end{equation}
	
\begin{corollary} \label{corlevset} 
 Under the hypotheses of Theorem \ref{teoest1}, suppose in addition that there exists $\lambda>0$ such that $\partial G_g(\lambda)\neq \emptyset$ and there exists $\gamma>0$ and $A>0$ such that if $|t-c|\leq \gamma$ then $d_H(\{g=c\},\{g=t\})\leq A|t-c|$. Then $d_H\big(\partial G_g(\lambda),\partial G_{\hat{g}_n}(\lambda)\big)=o(1/\beta_n) \quad \text{a.s.}
$
\end{corollary}

\begin{remark} As pointed out in \cite{cgmrc:06} section 2.4 point 1, the hypotheses of corollary \ref{corlevset} are fulfilled if $g$ is $C^2$ on a neighborhood $E$ of the level set $\lambda$ and the gradient of $g$ is strictly positive on $E$.
\end{remark}

\subsection{Level set estimation under shape restrictions}

\begin{figure}[ht]
\begin{center}
\includegraphics[height=3cm,width=3.8cm]{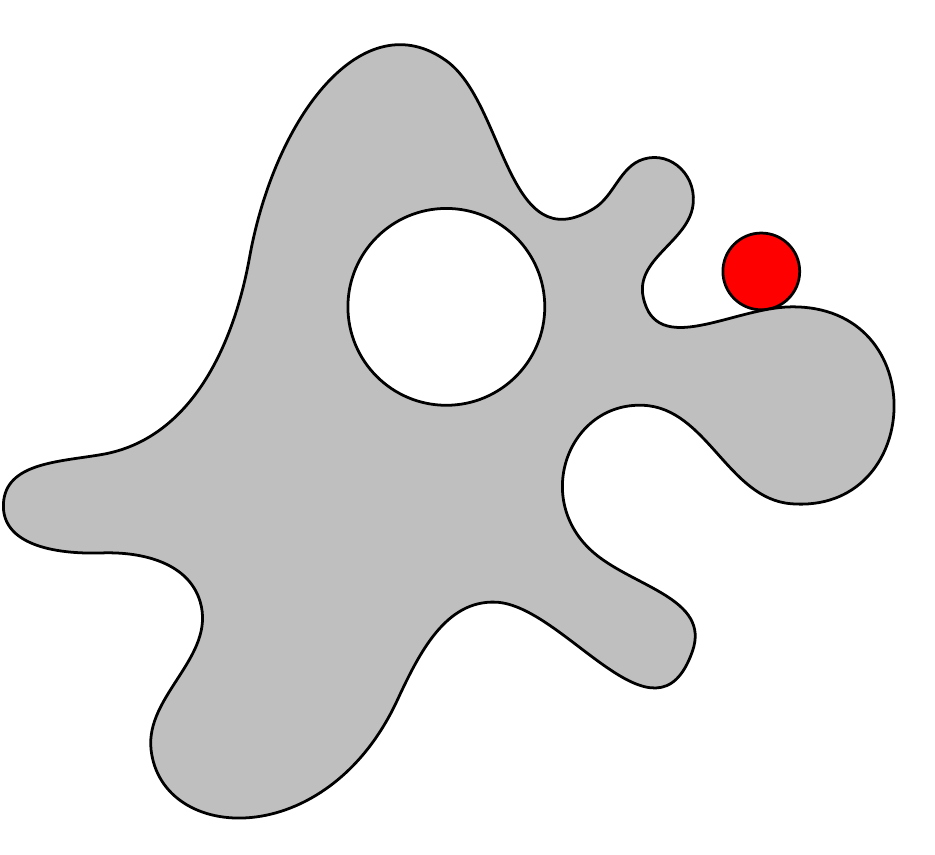}\\
\caption{A general $r$-convex set.
The small ball has radius $r$}\label{rconvex2}	
\end{center}
\end{figure}

In this subsection we propose another estimator of the level sets, 
under a quite general shape condition.
We assume that there exists an $r>0$ such that $\overline{G_g(\lambda)}$ is compact and $r$-convex, i.e. $\overline{G_g(\lambda)}=C_r\big(\overline{G_g(\lambda)}\big),$ where
$$
C_r\big(\overline{G_g(\lambda)}\big)=\bigcap_{\big\{ \mathring{\mathcal{B}}(x,r)\colon \mathring{\mathcal{B}}(x,r)\cap \overline{G_g(\lambda)}=\emptyset\big\}} \Big(\mathring{\mathcal{B}}(x,r)\Big)^c
$$
	is the $r$-convex hull of $G_g(\lambda)$.\\
This condition has been extensively studied in set estimation, see for instance \cite{cuevas:12}, \cite{pateiro:09} and \cite{rodriguez:07}. It is also related to the level set estimation problem, see \cite{walther:97}. Although $r$-convexity is much less restrictive than convexity, inlets that are too sharp are not allowed, see Figure \ref{rconvex2}.

Following the notation in \cite{federer:59}, let $\text{Unp}(S)$ be the set of points $x\in \mathbb{R}^d$ with a unique projection on $S$, denoted by $\xi_S(x)$. That is, for $x\in \text{Unp}(S)$, $\xi_S(x)$ is the unique point that attains the minimum of $\|x-y\|$ for $y\in S$. We write  $\delta_S(x)=\inf\{\|x-y\|:y\in S\}$.
\begin{definition} For $x\in S$, let {\it{reach}}$(S,x)=\sup\{r>0:\mathring{\mathcal{B}}(x,r)\subset\text{Unp}(S)\big\}$. The reach of $S$ is defined by $reach(S)=\inf\big\{reach(S,x):x\in S\big\},$ and $S$ is said to be of positive reach if $reach(S)>0$.
\end{definition}

The relation between $r$-convexity, reach, and rolling type conditions, has been studied in \cite{cuevas:12}.

\begin{definition}  The outer Minkowski content of $S\subset \mathbb{R}^d$  is given by 
	$$L_0(\partial S)=\lim_{\epsilon\rightarrow 0} \frac{\mu_L\big(B(S,\epsilon)\setminus S\big)}{\epsilon},$$
	provided that the limit exists and is finite.
\end{definition}

\begin{definition} Let $S\subset \mathbb{R}^d$ be a closed set. A ball of radius $r$ is said to roll freely in $S$ if for each boundary point $s\in \partial S$ there exists some $x\in S$ such that $s\in \mathcal{B}(x,r)\subset S$. The set $S$ is said to satisfy the outside $r$-rolling condition if a ball of radius $r$ rolls freely in $\overline{S^c}$.
\end{definition}

We will also assume the following condition.

\textbf{HR:} A level set $G_g(\lambda)$ fulfills \textbf{HR} if there exists $\delta_0>0$ and $r>0$ such that $\overline{G_g(\lambda+\eps)}$ is $r$-convex for all $-\delta_0<\eps<\delta_0$.

{Theorem 2} in \cite{walther:97} gives sufficient conditions for \textbf{HR} to hold, 
expressed in terms of the gradient of $g$.	{More precisely, it is shown the following result.
\begin{theorem} \label{thw97}
Let $g\colon \mathbb R^d \to \mathbb R$ and $-\infty< l \leq u < \sup g$. 
Assume that $g\in C^1(U)$ where $U$ is a bounded open set that contains $\overline{G_g(l-\eta)} \setminus G_g(u+ \eta)$ 
for some $\eta>0$; $\nabla g$ satisfies  $\Vert \nabla g\Vert \geq m >0$ on $U$ as well as a Lipschitz condition on $U$ (or on $\partial G_g(\lambda))$:  for all $\lambda \in (l,u)$ $\Vert \nabla g(x)-\nabla g(y)\Vert \leq k \Vert x-y \Vert,$  for $x, y \in U$ (or in $\partial G_g(\lambda))$. Then, for each $\lambda \in (l,u)$,  $\overline{G_g(\lambda)}$ and $G_g(\lambda)^c$ are $r_0$-convex with $r_0=m/k$.  
\end{theorem}}

\begin{lemma} \label{lemaux} 
Let $g\colon S\rightarrow \mathbb{R}$, where $S\subset \mathbb{R}^d$ is a compact set. 
Assume that $g\in C^2(S)$ and that $\lambda $ is such that there exists $0<\delta_1<\lambda$ for which
	$\nabla g(x) \neq 0$ for all $x\in \overline{G_g(\lambda - \delta_1)} \setminus G_g(\lambda+\delta_1)$. 
Then, 
for all $\eps<\delta_1$,
\begin{equation}\label{lemauxeq}
d_H\big(G_g(\lambda-\eps), G_g(\lambda+\eps)\big)\leq \frac{3M}{m^2}\eps,
\end{equation}
where 
$M=\max_{\{x\in \overline{G_g(\lambda - \delta_1)} \setminus G_g(\lambda+\delta_1)\}}\|\nabla g(x)\|$, 
and $m=\min_{\{x\in \overline{G_g(\lambda - \delta_1)} \setminus G_g(\lambda+\delta_1)\}}\|\nabla  g(x)\|$.
\end{lemma}

\begin{figure}[ht]
	\begin{center}
		\includegraphics[height=3cm,width=4cm]{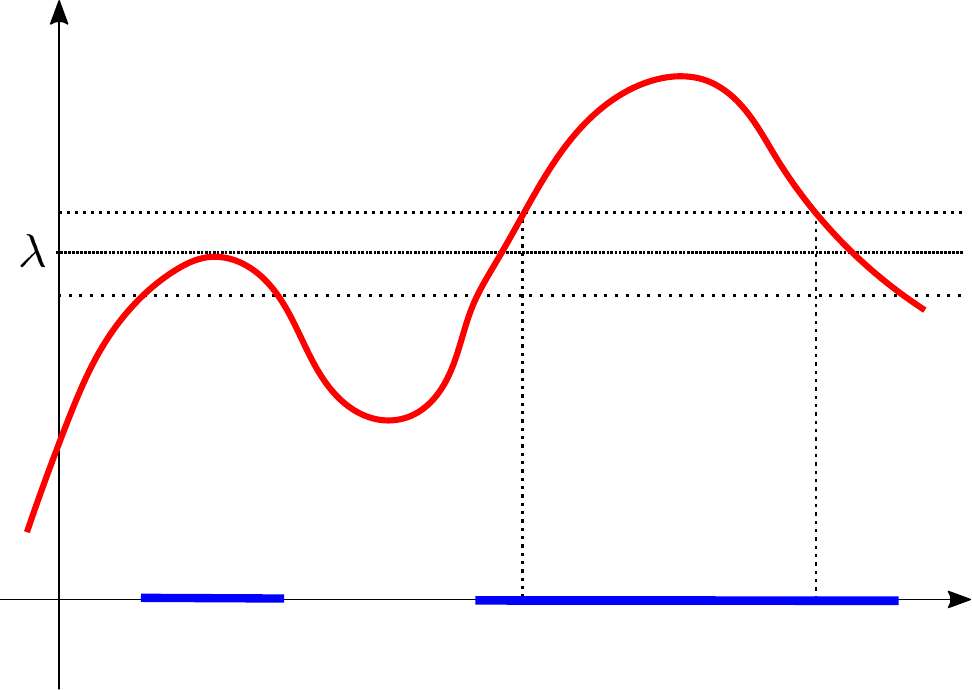}
	\end{center}
	\caption{ If $g'(x)=0$ for $x\in G_{g}(\lambda)$, not necessarily $d_H(G_g(\lambda+\eps),G_g(\lambda-\eps))\rightarrow 0$ as it is shown in the figure.}
	\label{HR}	
\end{figure}
Consider $\hat{g}_n$ as before. Assume that $g$ fulfills \textbf{HR}. 
We  study the convergence in the Hausdorff distance of the following estimator:
\begin{equation} \label{chull}
A_n(\lambda)=C_r\big(\{X_i\colon\hat{g}_n(X_i)>\lambda \}\big),
\end{equation}
i.e., the $r$-convex hull of the sample points beloging to the $\lambda$ level set of $\hat{g}_n$. 
The rates of convergence for the estimator \eqref{chull} in the independent case
were obtained in \cite{snrc:14}, where an estimator of the parameter $r$ was included. 
Observe that in our case it is not necessary to compute the whole set $G_{\hat{g}_n}(\lambda)$ 
(which in practice is not feasible in most cases), 
as the estimator proposed in Corollary \ref{corlevset} is based just on the sample points which belong to the set $G_{\hat{g}_n}(\lambda)$. 
Moreover, for the two dimensional case, the $r$-convex hull can be easily computed using the R software package {\tt alphahull} (see \cite{pateiro:10}).

\begin{theorem} \label{teoest2} 
Under the hypothesis of Theorem \ref{teoest1}, assume also that $g$ and $\lambda$ are in the hypothesis of Lemma \ref{lemaux}, 
and that condition \textbf{HR} holds; assume also that $0<g_0<g(x)$ for all $x\in S$. Let us denote $\Gamma_n$ the right hand side of \eqref{unifconvtasa}, with $\epsilon=1$. Let 
$\eps_n\rightarrow 0$ such that $\eps_n\beta_n>1$ for all $n$, assume also that $\eps_n<min\{\delta_0,\delta_1\}$  for all $n$, being $\delta_0$ as in condition \textbf{HR} and $\delta_1$ as in Lemma \ref{lemaux},  
Then, for all $n>n_2$ (where $n_2$ will be given in the proof), 
$$
\mathbb{P}\Big(d_H\Big(\overline{A_n(\lambda)},\overline{G_g(\lambda)}\Big)\leq \frac{3M}{m^2} \eps_n\Big)>1-3\Gamma_n.
$$
 \end{theorem}

The following Corollary follows directly from condition \textbf{HR} together with Theorem 3 in \cite{cuevas:12}.

\begin{corollary} Under the hypotheses of Theorem \ref{teoest2}, with probability one, 
	$$\lim_{n\rightarrow \infty} d_H(\partial \overline{A_n(\lambda)},\partial \overline{G_g(\lambda)})= 0.$$
\end{corollary}

\subsection{Estimation of level sets with a fixed content}

\begin{theorem} \label{thgtau}
Let $S\subset\mathbb{R}^d$ be a compact set and $\{X_n\}_{n\in \mathbb{N}}$ a geometrically ergodic Markov chain with state space $S$. 
For $\tau\in (0,1)$, define $l_\tau=\inf\{\lambda>0: \pi(G_g(\lambda))\leq 1-\tau\}$, $\pi$ being the stationary distribution. Assume that $\pi$ has a $C^2$ density $g$ such that $\|\nabla g(x)\|\neq 0$ for all $x\in U$, where $U$ is an open set containing $\overline{G_{g}(l_{\tau}-\eps_0)}\setminus G_g(l_\tau+\eps_0)$ for some $\tau>0$ and $0<\eps_0<l_\tau$. Let
$\hat{g}_n(x)=\frac{1}{n}\sum_{i=1}^nK_h(x-X_i)$
with $K$ a bounded Lipschitz density. Let $h=h_n$ be such that $h=\mathcal{O}((1/n)^{1/(d+1)})$. If we define 
$$\hat{l}_\tau=\inf\Big\{\lambda>0: \frac{1}{n}\#\big\{i:X_i \in G_{\hat{g}_n}(\lambda)\big\}\leq 1-\tau\Big\},$$
then, with probability one, $d_H\big(G_{\hat{g}_n}(\hat{l}_\tau),G_g(l_\tau)\big)\rightarrow 0.$
\end{theorem}

\subsection{Drift estimation}

{In what follows we propose an estimator of the drift function.
Assume that  $\{X_{t}:t\geq 0\}$ is uniformly sampled at times $\{t=t_1,t_2,\dots,t_n\}$ in the interval $[0,T]$, where $T>0$, i.e., a sample of size $n$ of the process $X_t$,  $\{X_{\Delta_{n,T}},X_{2\Delta_{n}},\dots,X_{n\Delta_{n,T}}\}$ is observed  at $t_i=i\Delta_{n,T}$ and $\Delta_{n,T}=T/n$. To simplify the notation we denote $\Delta$ instead of $\Delta_{n,T}$. 
 Let us fix $x\in int(S)$.	Denote by $N_x=\#\{1\leq i\leq n:X_{t_i}\in B(x,h_{n})\}$, for some $h_{n}\rightarrow 0$. 
 We define the estimator, 	
 \begin{equation*}
 \hat{\mu}_{n,T}(x)=\frac{1}{\Delta N_x} \sum_{i=1}^n (X_{t_{i+1}}-X_{t_i})\mathbb{I}_{\{X_{t_i}\in B(x,h_{n})\}}.
 \end{equation*} 
 \begin{theorem}\label{condrift}
 Assume that $T\rightarrow \infty$, $\Delta\rightarrow 0$, $h_{n}\rightarrow 0$, $\Delta nh_n^2\rightarrow \infty$. 
 Then, for all $x\in int(S)$
 \begin{equation}\label{consdrift}
 \hat{\mu}_{n,T}(x)\rightarrow \mu(x) \quad \text{ in probability.}
 \end{equation}
 \end{theorem}
 The proof is given in Appendix C.
According to Remark \ref{green}, in the gradient case,
the drift estimator can be  easily derived from the stationary density estimator, by using the plug-in rule 
\begin{equation}\label{driftest2}
\hat{\mu}_1(x)=\frac{1}{2} \nabla \log (\hat{g}_{n}(x)).
\end{equation}

 

\section{Examples}\label{examples}

In this section we first assess through a simulation study, the performance of the $r$-convex hull of the sample points belonging to the level set of the estimator, proposed in \eqref{chull}. 
Then we show the results of applying this method to real data.

\subsection{Simulations} \label{simulation}

The discrete version of the RBMD  \eqref{sde} is produced using the Euler scheme proposed in
\cite{bgt:04}, in the following way. We first choose a step $\delta>0$,  
and denote by $\operatorname{sym}(z)$ the symmetric
of the point $z$ with respect to $\partial S$.
We start with $X_0=x$ and suppose that we have obtained $X_i\in S$.  
To produce the following point, set
$$
Y_{i+1}=X_i+Z_i+\delta\mu(X_i),
$$ 
 where $Z_i$ is a centred Gaussian random vector,  independent w.r.t. $Z_1,\dots,Z_{i-1}$, 
with covariance matrix $\delta(I_d)_{\mathbb{R}^2}$.
Then
\begin{enumerate}
\item If $Y_{i+1}\in S$, set $X_{i+1}=Y_{i+1}$.
\item If $Y_{i+1}\notin S$ and $\operatorname{sym}(Y_{i+1})\in S$, set $X_{i+1}=\operatorname{sym}(Y_{i+1})$.
\item If $Y_{i+1}\notin S$ and $\operatorname{sym}(Y_{i+1})\notin S$, set $X_{i+1}=X_i$.
\end{enumerate}

In our example, we consider an RBMD in the set 
$S=E\setminus \mathcal{B}((4/5,0),1/2)$, where 
$E=\{(x,y)\in \mathbb{R}^2\colon 4x^2/9+y^2\leq 1\}$, with drift function given by  
$\mu(x,y)=-(x,y)$.
The stationary density is
\begin{equation}\label{case1f}
g(x)=\frac{1}{c}\exp\left[-(x^2+y^2)\right]\mathbb{I}_S(x,y)\quad \text{where } c=\iint_S \exp\left[-(x^2+y^2)\right]dxdy.
\end{equation}
The trajectory is shown in Figure \ref{c1fig1} for $\delta=0.001$ in the first row, and $\delta=0.003$ in the second row. The values for $N$ are $10,000; 50,000$ and $100,000$ in the first, second and third columns, respectively. 
 	
 	\begin{figure}[ht]
 		\begin{center}
 			\includegraphics[height=8.5cm,width=13.5cm]{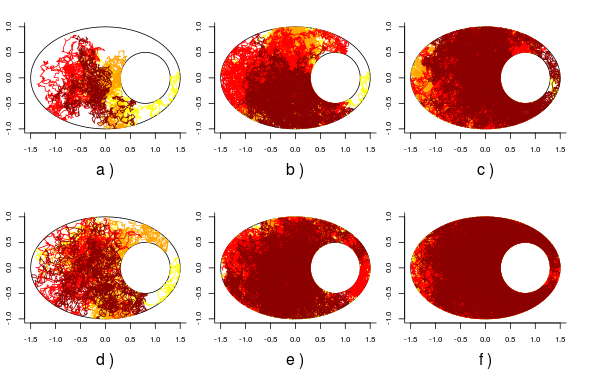}
 		\end{center}
 		\caption{\small The trajectory of the RBMD, for different values of $\delta$ and $N$, in a), b) and c) $\delta=0.001$ and $N=10,000$, $N=50,000$ and $N=100,000$, respectively. In d), e) and f), $\delta=0.003$ and $N=10,000$, $N=50,000$ and $N=100,000$, respectively.}
 		\label{c1fig1}
 	\end{figure}
 	The function \eqref{case1f} is shown in Figure \ref{c1fig2} a), while in b) there is shown the estimated density using a Gaussian kernel with bandwidth $h=0_\cdot 2$; in c) there is shown the estimated density using an Epanechnikov kernel with bandwidth $h=0_\cdot 4$. 
 		 In both cases we have used the trajectory shown in Figure \ref{c1fig1}, with $\delta=0.003$ and $N=100,000$.  Since we can estimate the support, we have forced the estimation to be 0 outside the estimation of the support. 
 	\begin{figure}[ht]
 		\begin{center}
			\includegraphics[width=.9\textwidth]{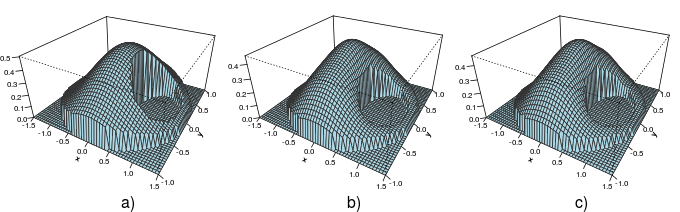}
 		\end{center}
 		\caption{a) Real density, b) estimated using Gaussian Kernel with  $h=0_\cdot 2$, c) estimated using Epanechnikov kernel with $h=0_\cdot 4$.}
 		\label{c1fig2}
 	\end{figure}
 	 
For the level sets, we have considered the levels $\lambda =$ 0.44, 0.41, 0.34, 0.27 and 0.03. Figure \ref{c1fig4} a) shows the theoretical level sets for the considered values of $\lambda$, while in b) there are shown the corresponding estimated level sets. The estimation is based on the trajectory with $\delta=0.003$ and $N=500,000$ using \eqref{chull} with $r=0.4$. We have used the Gaussian kernel with $h=0_\cdot 1$. The choice of an optimal bandwidth for level set estimation has been studied recently for the iid case, see \cite{q:18}. Although it should behaves similarly for geometric mixing processes, to extend the results in \cite{q:18} is far beyond the aim of this paper. It is clear that the hole in the domain will produce border effects for the density estimation, and therefore for the level sets. A way to overcome this problem (which is computationally very expensive) is to first estimate the support using the $r$-convex hull of the trajectory and then use a variable bandwidth kernel estimate where the bandwidth is given by the lesser of a fixed $h$ and the distance from the point $x$ to the boundary of the support. 
	
 	\begin{figure}[h]
 		\begin{center}
 			\includegraphics[height=5cm,width=6.5cm]{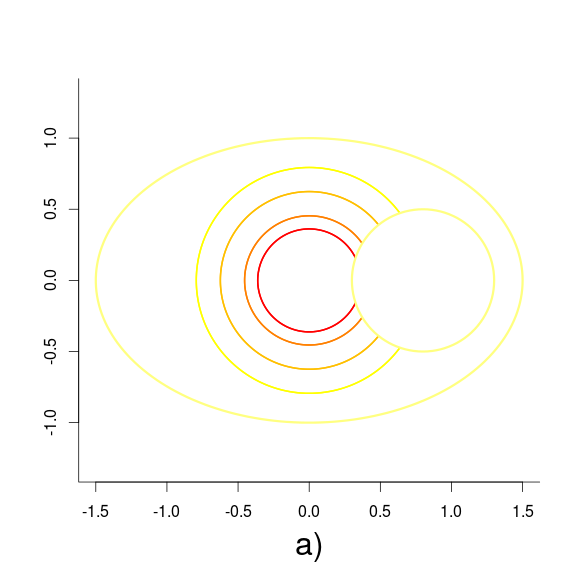}			 												\includegraphics[height=5cm,width=6.5cm]{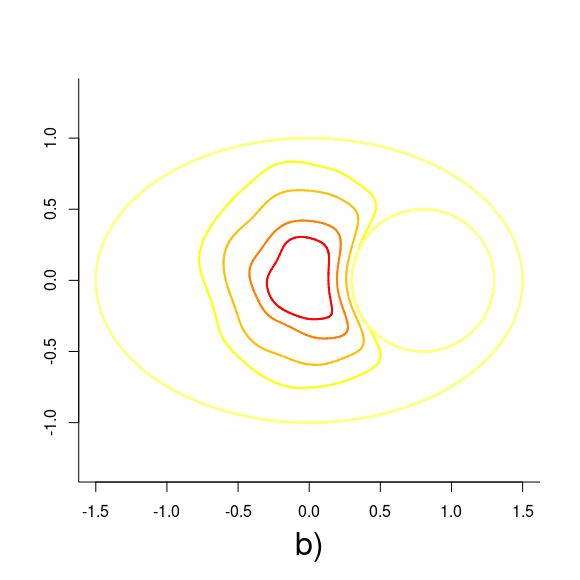}						
 			\end{center}
 		\caption{a) Theoretical level sets. b) Estimation using \eqref{chull} for $r=0.4$, with Gaussian kernel and a bandwidth $h=0.1$. In red the core-area.}
 		\label{c1fig4}
 	\end{figure}

In Figure \ref{vf1} a) we represent the theoretical vector field corresponding to the drift, while in b)  we provide the estimator \eqref{driftest2} based on the trajectory given in Figure \ref{c1fig1} f), using the Gaussian kernel and a bandwidth $h=0.45$.
\begin{figure}[h]
 		\begin{center}
 			\includegraphics[height=4.5cm,width=6.6cm]{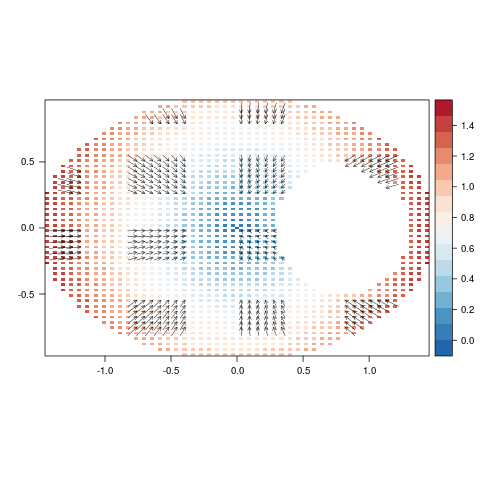} \hfill
 			\includegraphics[height=4.5cm,width=6.6cm]{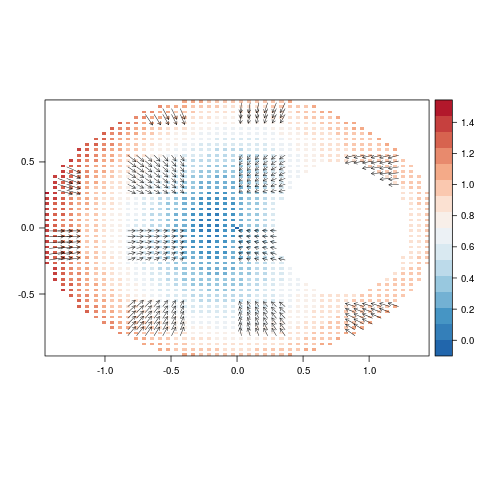}
 		\end{center}
 		\caption{{(Left) Theoretical drift  (Right) estimation with \eqref{driftest2}}}
 		\label{vf1}
 	\end{figure}

\subsection{Real data examples} \label{realdata}

We considered a dataset from the Movebank database, where a natural
barrier acts as a boundary of the animal's movement. GPS collars were placed on elephants
in Loango National Park in western Gabon. The area is protected by the Atlantic Ocean
on the west and by Lagoon Igu\'ela on the east. Figure \ref{eleph} a) shows in red the movement of an elephant
with estimator $N = 1633$ for recorded positions. In blue we represent the boundary of the $r$-convex hull
estimator for $r=0.02$. The estimated density is shown in b), using the Gaussian kernel with bandwidth $h=0.01$. 
The $r$-convex hulls of the level sets are shown in c) for $\lambda_1=100,\lambda_2=600$, $\lambda_3=1100$, $\lambda_4=1600$,and $r=0.02$.
In d) we represent the estimation of the drift, using \eqref{driftest2} with $h=0.5$.
 \begin{figure}[ht]
	\begin{center}
 		\includegraphics[height=4cm,width=5cm]{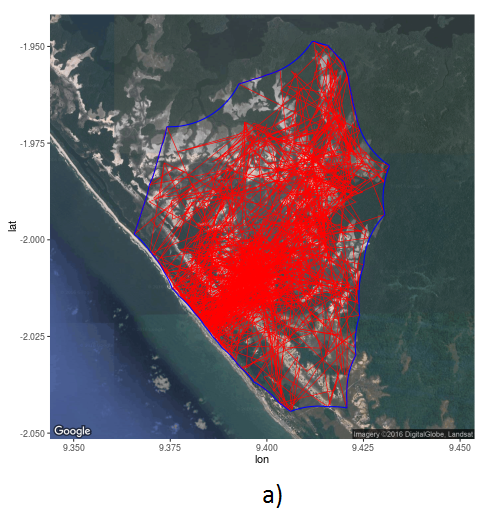}
 	    \includegraphics[height=4cm,width=5cm]{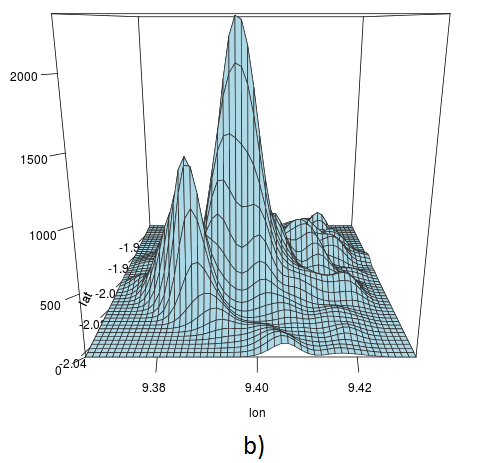}\\
 	 \includegraphics[height=4cm,width=5cm]{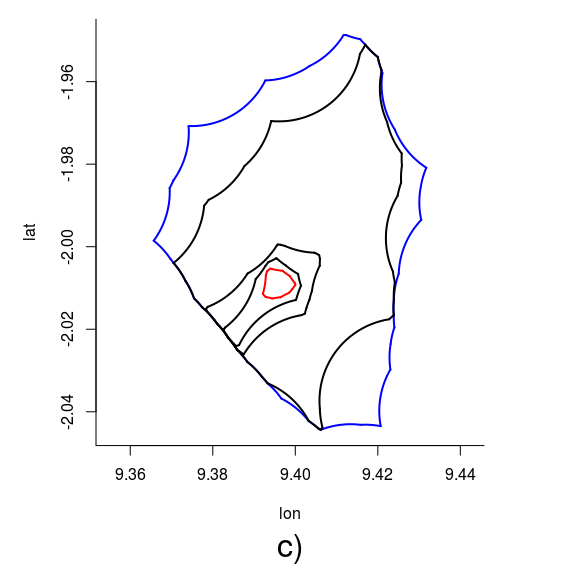}
  		\includegraphics[height=4cm,width=5cm]{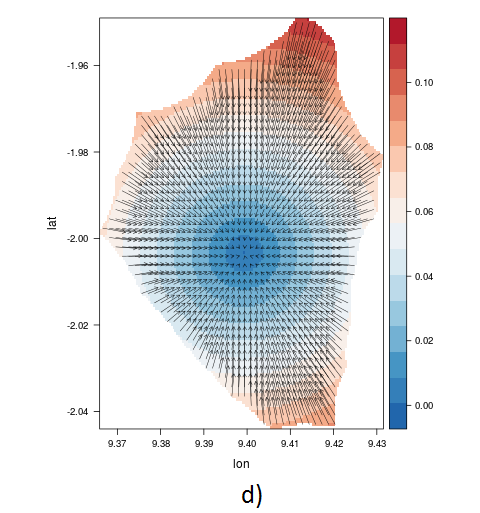}
 		\end{center}
 {\caption{\small  a) Trajectory and $0.02$-convex hull b) density estimator using Gaussian Kernel with $h=0.01$ c) $r$-convex hull of  level sets d) Estimation of the drift}}
 		\label{eleph}
 	\end{figure} 		

 \section{Appendix A}
 
Here we include the proofs of the propositions stated in Section \ref{driftest}.

\textit{Proof of Theorem \ref{teoest1}}

We will deal separately with each term on the right hand side of the following inequality:
	\begin{equation} \label{unifeq1}
	\beta_n\sup_x \big|g(x)-\hat{g}_n(x)\big|\leq \beta_n\sup_x \big|g(x)-\EXP_\pi(\hat{g}_n(x))\big|+\beta_n\sup_x \big|\hat{g}_n(x)-\EXP_\pi(\hat{g}_n(x))\big|.
	\end{equation}
	First we bound the bias term.  Let $C_g$ such that $|g(x)-g(y)|\leq C_g\|x-y\|$. Then
	\begin{multline} \label{eq1}
	\big|\EXP_\pi(K_h(x-X_k))-g(x)\big|\leq \left|\int_S K_h(x-y)g(y)dy-g(x)\right|\\ \leq 
	\int_{S} K_h(x-y)|g(y)-g(x)|dy\leq hC_g\int_{\mathbb{R}^d} \|u\|K(u)du=\kappa C_g h.
	\end{multline}
	Observe that $\EXP_{\mu_0}(K_h(x-X_k))=\int_S K_h(x-y)\mu_k(dy)$. 
	Recall that $k_1=\max_x K(x)$. Now, by \eqref{ge},
	\begin{equation} \label{eqth2}
	\left|\int_S K_h(x-y)\mu_k(dy)-\int_S K_h(x-y)d\pi(y)\right|\leq\|K_h(x-y)\|_\infty \|\mu_k-\pi\|_{TV} \leq \frac{k_1}{h^d}\gamma\rho^k.
	\end{equation}
	Observe that $\EXP(\hat{g}_n(x))=\frac{1}{n}\sum_{k=1}^n \EXP_{\mu_0}(K_h(x-X_k))$, and $\EXP_\pi (K_h(x-X_k))=\int_S K_h(x-y)g(y)dy$. Hence \eqref{eqth2} implies
	\begin{multline} \label{eqth2.1}
	\beta_n\sup_x\left|\frac{1}{n}\sum_{k=1}^n\big[\EXP_{\mu_0}(K_h(x-X_k))-\EXP_\pi(K_h(x-X_k))\big]\right|\leq \frac{\beta_n}{nh^d} 2k_1\gamma  \sum_{k=1}^n\rho^k\leq \\
		\frac{\beta_n}{nh^d}  2k_1\gamma  \sum_{k=1}^\infty\rho^k=C\frac{\beta_n}{nh^d}		\end{multline}
	which, together with \eqref{eq1}, implies that 
	\begin{equation} \label{cotasesgo}
	\beta_n\sup_x \big|\EXP(\hat{g}_n(x))-g(x)\big|\leq C \frac{\beta_n}{nh^d}+\kappa C_g h\beta_n.
	\end{equation}

	It remains to prove that $\beta_n\sup_x \big|\hat{g}_n(x)-\EXP(\hat{g}_n(x))\big|\rightarrow 0$.
	Since $S$ is compact, we can cover $S$ with $\nu\leq \frac{c}{h^{d(d+2)}}$ balls of radius $h^{d+2}$ centred at some fixed points $\{x_1,\dots,x_\nu\}\subset S$, $c$ being a positive constant depending only on $d$ and $\mu_L(S)$. For $i=1,\dots,\nu$,
	\begin{equation*} \mathbb{P}\Big(\big|\hat{g}_n(x_i)-\EXP(\hat{g}_n(x_i))\big|>\eps\Big)=\mathbb{P}\left(\left|\sum_{j=1}^n\big[K_h(x_i-X_j)-\EXP_{\mu_0}(K_h(x_i-X_j))\big]\right|>n\eps\right).
	\end{equation*}
	By proposition 4.1 of \cite{cd:05}, the sequence $\{X_n\}_{n\in\mathbb{N}}$ is $\varphi$ mixing with 
	$\varphi(n)=2\gamma \rho^n$ ($\gamma$ as in \eqref{ge}). 
	Let $x\in S$ and $x_i$ be such that $\|x-x_i\|<h^{d+2}$.  
	Then, since $K$ is Lipschitz, denote by $R$ the Lipschitz constant of $K$, then
	$$
	|\hat{g}_n(x)-\hat{g}_n(x_i)|\leq \frac{1}{nh^d}\sum_{j=1}^n \left|K\left(\frac{x-X_j}{h}\right)-K\left(\frac{x_i-X_j}{h}\right)\right|\leq \frac{1}{h^{d+1}}R\|x-x_i\|\leq Rh.
	$$
	Hence, $\big|\hat{g}_n(x)-\EXP(\hat{g}_n(x))\big|\leq \big|\hat{g}_n(x_i)-\EXP(\hat{g}_n(x_i))\big|+2Rh.$
	If we take $n$ so large that $2Rh<\eps/(2\beta_n)$, we get 
	\begin{equation*}
	\mathbb{P}\left(\sup_{x\in \mathcal{B}(x_i,h^{d+2})}\big|\hat{g}_n(x)-\EXP(\hat{g}_n(x))\big|>\frac{\eps}{\beta_n}\right)\leq \mathbb{P}\left(\left|\hat{g}_n(x_i)-\EXP(\hat{g}_n(x_i))\right|>\frac{\eps}{2\beta_n}\right).
	\end{equation*}
	
	Now use the Bernstein inequality \eqref{bern} with $Y_j=K((x-X_j)/h)-\EXP_{\mu_0}(K(x-X_j)/h)$ and $C_1=2k_1$. Recall that $g_1=\max_{x\in S}g(x)$. Let us take $n_0$ such that for all $n>n_0$ $\rho^n/(h_n^d g_1\mu_L(S))<1$ and $2Rh<\eps/(2\beta_n)$. Denote by $C''=2k_1\gamma g_1\mu_L(S)$, then for $n>n_0$, by \eqref{eqth2}, 
	$$\mathbb{E}_{\mu_0}\left(K\left(\frac{x-X_j}{h}\right)\right)\leq k_1\gamma \rho^n+ \int_SK\left(\frac{x-y}{h}\right)g(y)dy\leq k_1\gamma \rho^n+k_1h^dg_1\mu_L(S)\leq C''h^d.$$
	Hence $\eta=2C''h^{d}$, $D\leq 2k_1C'' h^{d}$ , and $\tilde{\varphi}(m)\leq \sum_{i=1}^\infty 2\gamma \rho^i= 2\gamma \rho/(1-\rho)<2\gamma$,  so 
	$C_2=12k_1C''h^d(1+16\gamma)$. Since $\alpha_n=o(1/\beta_n)$, if $m=\lfloor \beta_n\rfloor$, then $\alpha_n C_2 h^{-d} < \varepsilon/(4\beta_n)$ and $\alpha_n m C_1< 1/4$ for $n$ large enough. On the other hand,  since $\log(n)/\beta_n\to 0$,
	$$
	3e^{1/2}n\frac{\varphi(m)}{m} \to 0, \text{ as }n\rightarrow \infty.
	$$
	Let us take $n_1>n_0$, such that for all $n>n_1$, $\beta_n\alpha_n C_2 h^{-d}<\epsilon/4$, $\alpha_n m C_1< 1/4$ and $2\exp\Big(3e^{1/2}n\frac{\varphi(m)}{m}\Big)<3$.
	
	Now the Bernstein inequality  implies that, for all $n>n_1$	
	\begin{align*} 
	\mathbb{P}\left(\big|\hat{g}_n(x_i)-\EXP(\hat{g}_n(x_i))\big|>\frac{\eps}{2\beta_n}\right)&\leq 2\exp\Big(3e^{1/2}n\frac{\varphi(m)}{m}\Big)\exp\Big(-\frac{\alpha_n\eps n h^d}{2\beta_n} +\alpha_n^2C_2n\Big)\nonumber\\
	&\leq 3 \exp\Big(-\frac{\eps\alpha_n n h^{d}}{4\beta_n}\Big).
	\end{align*}
	Finally, for $n>n_1$,
	\begin{align*} 
	\mathbb{P}\Big(\sup_x\big|\hat{g}_n(x)-\EXP(\hat{g}_n(x))\big|>\frac{\eps}{\beta_n}\Big)
	&\leq \sum_{i=1}^{\nu} \mathbb{P}\left(\big|\hat{g}_n(x_i)-\EXP(\hat{g}_n(x_i))\big|>\frac{\eps}{2\beta_n}\right)\nonumber \\
	&\leq  3\frac{c}{h^{d(d+2)}}\exp\Big(- \frac{\epsilon \alpha_n n h^{d}}{4\beta_n}\Big)\\
	&\leq  3c\exp\Big(-\epsilon \alpha_n\frac{nh^d}{4\beta_n}-(d(d+2))\log(h)\Big),
	\end{align*}
which, together with \eqref{unifeq1}  and \eqref{cotasesgo}, implies \eqref{unifconvtasa}. \\
To prove the almost surely convergence, just observe that $\alpha_nnh^d/(4\beta_n \log(n))\rightarrow \infty$ and $nh\rightarrow \infty$, imply that
$$\frac{1}{\log(n)}\Big[\epsilon\alpha_n\frac{nh^d}{4\beta_n}+(d(d+2))\log(h)\Big]\rightarrow \infty,$$
and then we can apply Borel-Cantelli Lemma.

\textit{Proof of Lemma \ref{lemaux}}

	Let $x\in G_{g}(\lambda-\eps)$, $y_t=x+t\nabla g(x)$ and $t=3\eps /m^2$.
	We have $\|y_t-x\|<\frac{3}{m^2}\eps M$. 
	To prove \eqref{lemauxeq} it is enough to verify that $y_t\in G_g(\lambda+\eps)$. 
	From a Taylor expansion at $x$,  we obtain that for some $\theta \in [x,y_t]$:
	$$
	\aligned
	g(y_t)	&=g(x)+\nabla g(x)^T(y_t-x)+\frac12(y_t-x)^T H_\theta (y_t-x)\\
	&> \lambda-\eps+\frac{3\eps}{m^2}\|\nabla g(x)\|^2+\frac{9\eps^2}{2m^4}\nabla g(x)^T H_\theta \nabla g(x),
	\endaligned
	$$
	where $H_\theta$ is the Hessian matrix of $g$ at $\theta$. 
	Since $g$ is $C^2$, there exists a constant $C>0$ such that $|\nabla g(x)^T H_\theta \nabla g(x)|\leq C\|\nabla g(x)\|^2$, 
	from where it follows that for $\eps<2m^4/(9M^2C)$,
	$$g(y_t)> \lambda+2\eps-\frac{9M^2C}{2m^4}\eps^2\geq \lambda +\eps,
	$$	 
	and $y_t\in G_g(\lambda+\eps)$, concluding the proof.

\textit{Proof of Theorem \ref{teoest2}}

	Let us consider $n>n_1$ where $n_1$ is given in Theorem \ref{teoest1}. Let us denote $\mathcal{A}_n=\{ \beta_n\sup_{x}|\hat{g}_n(x)-g(x)|<1\}$, we now that $\mathbb{P}(\mathcal{A}_n)>1-\Gamma_n$ for all $n>n_1$. Since $\eps_n\beta_n>1$ and $\eps_n<\delta_0$, by condition \textbf{HR}, $\mathbb{P}(A_n(\lambda)\subset G_g(\lambda-\eps_n))>1-\Gamma_n$ for all  $n>n_1$. By Lemma \ref{lemaux} it is enough to prove that there exists $n_2$ such that for all $n>n_2$, $\mathbb{P}(G_g(\lambda+\eps_n)\subset A_n(\lambda))>1-2\Gamma_n$ or what is the same
	we have to prove that for $n>n_1$, $\mathbb{P}(\exists x_n \in G_g(\lambda+\eps_n):x_n\notin A_n(\lambda))\leq 2\Gamma_n$. \\
	Let us denote $\mathcal{X}_n=\{X_1,\ldots,X_n\}$
	$$\mathcal{C}_n=\Big\{\exists x_n \in G_g(\lambda+\eps_n) \text{ such that }\exists y_n: x_n\in \mathcal{B}(y_n,r),\#\big\{\mathcal{X}_n\cap  \mathcal{B}(y_n,r)\}=0\Big\}.$$
	Then,  
	\begin{multline*}
	\big\{\exists x_n \in G_g(\lambda+\eps_n):x_n\notin A_n(\lambda)\big\}\subset \mathcal{C}_n\cup\\ 
	\Big\{\{\exists x_n \in G_g(\lambda+\eps_n) \text{ such that }\exists y_n: x_n\in \mathcal{B}(y_n,r),\#\big\{\mathcal{X}_n\cap  \mathcal{B}(y_n,r)\}>0\big\}\\
	\cap \big\{\forall X_i\in \mathcal{B}(y_n,r), \hat{g}_n(X_i)\leq \lambda\big\}\Big\}=\mathcal{C}_n\cup \mathcal{F}_n.
	\end{multline*}
	Since $g$ is Lipschitz (denote by $C_g$ the Lipschitz constant) if $x_n\in G_g(\lambda+\eps_n)$, $g(z)>\lambda+\eps_n/2$ for all $z\in \mathcal{B}(x_n,\nu_n)$ where $\nu_n=\eps_n/(2C_g)$. Then on $\mathcal{A}_n$, for all $n>n_1$ $\hat{g}_n(z)>\lambda$, for all $z\in \mathcal{B}(x_n,\nu_n)$. Then 
	\begin{multline*}
	\mathcal{E}_n=\Big\{\exists x_n \in G_g(\lambda+\eps_n) \text{ such that }\exists y_n: x_n\in \mathcal{B}(y_n,r),\#\{ \mathcal{X}_n\cap  \mathcal{B}(y_n,r)\}>0\Big\}\\ \cap \big\{\forall X_i\in \mathcal{B}(y_n,r)\cap \mathcal{B}(x_n,\nu_n), \hat{g}_n(X_i)\leq \lambda\big\} \subset \mathcal{A}_n^c.
	\end{multline*}
	And then, $\mathbb{P}(\mathcal{F}_n)\leq \mathbb{P}(\mathcal{E}_n)\leq \Gamma_n$.
	
	 Let us bound $\mathbb{P}(\mathcal{C}_n)\leq \Gamma_n$.   To do that, let us introduce, for each fixed $n>n_1$,the random variables 
	$$Z_k(y)=K(\|X_k-y\|/r) \quad k=1,\dots,n,
	$$
	where $K$ is a Lipschitz function such that $\mathbb{I}_{[0,1/2]}(x)\leq K(x)\leq \mathbb{I}_{{[0,1]}}(x)$ and $K(x)>0$ for all $x\in (1/2,1)$, then
	$$\mathbb{P}(\mathcal{C}_n)\leq \mathbb{P}\Big( \inf_{\{y\in S\}} \frac{1}{n}\sum_{k=1}^n Z_k(y)= 0 \Big).$$
		Proceeding as in $\eqref{eqth2.1}$,
	\begin{equation} \label{bound0}
		\sup_{y\in S} \left|\frac{1}{n}\sum_{k=1}^n\big[\EXP_{\mu_0}(Z_k(y))-\EXP_\pi(Z_k(y))\big]\right|\leq \frac{2\gamma}{n} \sum_{k=1}^\infty \rho^k.
		\end{equation}
    Since $g(x)>g_0>0$ for all $z\in S$,
	\begin{equation} \label{boundex}
	\mathbb{E}_{\pi}(Z_k(y))=\int_{\mathcal{B}(y,r)}K(\|t-y\|/r)g(t)dt\geq g_0(r/2)^d\omega_d>0.
	\end{equation}
   Let us fix $0<\epsilon<g_0(r/2)^d\omega_d/3$, from  \eqref{bound0} and \eqref{boundex}, if we take $n$ large enough such that $\frac{2\gamma}{n} \sum_{k=1}^\infty \rho^k<g_0(r/2)^d\omega_d/3$, it is enough to prove that there exists $n_2$ such that for all $n>n_2$,
  	$$\mathbb{P}\Big( \sup_{y\in S} \Big|\frac{1}{n}\sum_{k=1}^n [Z_k(y)-\mathbb{E}_{\mu_0}(Z_k(y))]\Big|> \epsilon\Big)<\Gamma_n.$$
  	As before, since $S$ is compact, we can cover it with $\zeta\leq c/\iota^d$ balls of radius $\iota$ centred at some fixed points $\{x_1,\dots,x_{\zeta}\}$ where $c$ is a constant which depends only on $d$ and $\mu_L(S)$.
First, observe that if $\frac{2\gamma}{n} \sum_{k=1}^\infty \rho^k<\epsilon/5$, then, 
\begin{align*}
\sup_{y_i} \Big|\frac{1}{n}\sum_{k=1}^n [\mathbb{E}_{\mu_0}(Z_k(y_i))-\mathbb{E}_{\pi}(Z_k(y_i))]\Big|&\leq \frac{\epsilon}{5}.
\end{align*}  
And if $\mu_L(B(y_i,r)\triangle B(y,r))g_1<\epsilon/5$, where $g_1$ is the maximum of $g$,
\begin{align*}
\sup_{y_i\in B(y,\iota)} \Big|\frac{1}{n}\sum_{k=1}^n [\mathbb{E}_{\pi}(Z_k(y_i))-\mathbb{E}_{\pi}(Z_k(y))]\Big|&\leq \frac{\epsilon}{5}\\
\sup_{y\in S} \Big|\frac{1}{n}\sum_{k=1}^n [\mathbb{E}_{\pi}(Z_k(y))-\mathbb{E}_{\mu_0}(Z_k(y))]\Big|\leq \frac{\epsilon}{5}.
\end{align*}

Using Berstein inequality, as in Theorem 1, we can bound, for a fixed $y$,  
\begin{equation} \label{bern2}
\mathbb{P}\Big(\Big|\frac{1}{n}\sum_{k=1}^n [Z_k(y)-\mathbb{E}_{\mu_0}(Z_k(y))]\Big|>\frac{\epsilon}{\beta_n}\Big)\leq 3\exp\Big(-\frac{\epsilon\alpha_n n}{4\beta_n}\Big).
\end{equation}
where $\alpha_n=o(1/\beta_n)$ and $\log(n)/\beta_n\rightarrow 0$. 
Then
\begin{multline*}
\mathbb{P}\Big( \sup_{y\in S} \Big|\frac{1}{n}\sum_{k=1}^n [Z_k(y)-\mathbb{E}_{\mu_0}(Z_k(y))]\Big|>\epsilon \Big)\leq
\mathbb{P}\Big( \sup_{y_i\in B(y,\iota)} \Big|\frac{1}{n}\sum_{k=1}^n [Z_k(y)-Z_k(y_i)]\Big|>\frac{\epsilon}{5} \Big)+\\
\mathbb{P}\Big( \sup_{y_i} \Big|\frac{1}{n}\sum_{k=1}^n [Z_k(y_i)-\mathbb{E}_{\mu_0}(Z_k(y_i))]\Big|>\frac{\epsilon}{5} \Big)=I_1+I_2.
\end{multline*}

Since $K$ is Lipschitz  (let us denote $C_K$ the Lipschitz constant of $K$) we can bound $I_1\leq C_K\iota/r$ and from \eqref{bern2}, 
$I_2\leq \frac{3c}{\iota^d}\exp\Big(-\frac{\epsilon\alpha_n n}{4\beta_n}\Big)$. Now take $\iota=h_n$ (being $h_n$ as in Theorem 1), then
 $I_1+I_2\leq \Gamma_n$.

In order to prove Theorem  \ref{thgtau} we will need two lemmas. For the first, recall that given a probability distribution $P$,  $\mathcal{A}$ is a $P$-uniformity class if $\sup_{A\in \mathcal{A}}|P_n(A)-P(A)|\rightarrow 0$ whenever $P_n\rightarrow P$ weakly. 
Theorem 5 in \cite{cuevas:12} proves that the class of sets with reach bounded from below by a positive constant included in a compact set is a $P$-uniformity class.

\begin{lemma} \label{lempcont} Let $S\subset\mathbb{R}^d$ be a compact set and $g:S\rightarrow \mathbb{R}$ a $C^2$ function such that that there exists an $\eps_0>0$ and a $c>0$ such that $\|\nabla g(x)\|>m$ for all $x\in U$,  where $U$ is an open set containing $\overline{G_{g}(l_{\tau}-\eps_0)}\setminus G_g(l_\tau+\eps_0)$. Then $\{G_g(\lambda):l_\tau-\eps_0/2 \leq\lambda\leq l_\tau+\eps_0/2\}$ is a $P$-uniformity class for all probability distributions $P$ on $S$ absolutely continuous w.r.t. Lebesgue measure.
	\begin{proof} It is enough to prove that there exists an $r>0$ such that for all $l_\tau-\eps_0 <\lambda<l_\tau+\eps_0$, $\emph{reach}(G_g(\lambda))>r>0$. By Theorem \ref{thw97} and theorem 1 of \cite{walther:99}, there exists an $r>0$ such that for all $l_\tau-\eps_0 <\lambda<l_\tau+\eps_0$, $G_g(\lambda)$ satisfies the inner and outer $r$-rolling conditions. This together with lemma 2.3 in \cite{pateiro:09} implies that $\emph{reach}(G_g(\lambda))>r>0$ for all  $l_\tau-\eps_0/2\leq \lambda\leq l_\tau+\eps_0/2$. 
	\end{proof}
\end{lemma}

The following Lemma can be derived from Lemma 2b) in \cite{walther:97}, for the sake of completeness we keep the proof, which is a straightforward consequence of Lemma \ref{lemaux}. 

\begin{lemma} \label{lemdhaus} Under the hypotheses of Lemma \ref{lempcont}, for all $0\leq \eps<\eps_0/2$ and all $l_\tau-\eps<\lambda <l_\tau+\eps$,
	\begin{equation*}
	G_g(\lambda-\eps)\setminus G_g(\lambda+\eps)\subset B\Big(\partial G_g(\lambda),  \frac{3\eps M}{m^2}\Big),
	\end{equation*}
	where $M=\max_{\{x\in \overline{G_g(l_\tau- \eps_0)} \setminus G_g(l_\tau+\eps_0)\}}\|\nabla g(x)\|$ and $m=\min_{\{x\in \overline{G_g(\lambda - \delta_1)} \setminus G_g(\lambda+\delta_1)\}}\|\nabla  g(x)\|$.
	\begin{proof} By Lemma \ref{lemaux}, for all $\eps<\eps_0/2$ and all $l_\tau-\eps<\lambda <l_\tau+\eps$, 
		$d_H(G_g(\lambda+\eps),G_g(\lambda-\eps))\leq 3\eps M/m^2.$
		If we take $x\in G_g(\lambda-\eps)$ with $g(x)\leq \lambda$ and $y\in G_g(\lambda+\eps)$, then there exists a $t\in [x,y]$ (the segment joining $x$ and $y$) such that $g(t)=\lambda$, and so $t\in \partial G_g(\lambda)$, which concludes the proof.
	\end{proof}
\end{lemma}

\textit{Proof of Theorem \ref{thgtau}}\\
By Remark \ref{remconf} ii) we have that $\sup_{x\in S}|\hat{g}_n(x)-g(x)|\rightarrow 0$ a.s. We will prove that $\hat{l}_\tau\rightarrow l_\tau$ a.s. Define $L(\lambda)=\pi(G_g(\lambda))$, $\hat{L}(\lambda)=\frac{1}{n}\#\{i:X_i\in G_{\hat g_n}(\lambda)\}$ and $\tilde{L}(\lambda)=\frac{1}{n}\#\{i:X_i\in G_{g}(\lambda)\}$. Write
\begin{equation*}
\sup_{l_\tau- \frac{\eps_0}{2}\leq \lambda\leq l_\tau+\frac{\eps_0}{2}}|L(\lambda)-\hat{L}(\lambda)|\leq \sup_{l_\tau- \frac{\eps_0}{2}\leq \lambda\leq l_\tau+\frac{\eps_0}{2}}|L(\lambda)-\tilde{L}(\lambda)|+\sup_{l_\tau- \frac{\eps_0}{2}\leq \lambda\leq l_\tau+\frac{\eps_0}{2}}|\tilde{L}(\lambda)-\hat{L}(\lambda)|.
\end{equation*}
\begin{align*}
|\tilde{L}(\lambda)-\hat{L}(\lambda)|=&\frac{1}{n}\Big|\#\{i:X_i\in G_{g}(\lambda)\}-\#\{i:X_i\in G_{\hat g_n}(\lambda)\}\Big|\\
=&\frac{1}{n}\Big(\#\{i:X_i\in G_{g}(\lambda)\setminus G_{\hat{g}_n}(\lambda)\}+\#\{i:X_i\in G_{\hat{g}_n}(\lambda)\setminus G_{g}(\lambda)\}\Big).
\end{align*}
Since $\sup_x|\hat{g}_n(x)-g(x)|\rightarrow 0$ a.s., we have that for all $\lambda$ and $\eps$, $G_g(\lambda+\eps)\subset G_{\hat{g}_n}(\lambda)\subset G_g(\lambda-\eps)$ with probability one, for $n$ large enough. Then, with probability one, for $n$ large enough, for all $0\leq \eps<\eps_0/2$,
$$\sup_{l_\tau- \frac{\eps_0}{2}\leq \lambda \leq l_\tau+\frac{\eps_0}{2}}|\tilde{L}(\lambda)-\hat{L}(\lambda)|\leq\sup_{l_\tau- \frac{\eps_0}{2}\leq \lambda\leq l_\tau+\frac{\eps_0}{2}} \ \ \frac{2}{n}\#\Big\{i:X_i\in G_{g}(\lambda-\eps)\setminus G_g(\lambda+\eps)\Big\}.$$
By Lemma \ref{lempcont}, $G_g(\lambda)$ is a $P$-uniformity class. Hence,
$$\sup_{l_\tau- \frac{\eps_0}{2}\leq\lambda\leq l_\tau+\frac{\eps_0}{2}}\Big|\frac{1}{n}\#\Big\{i:X_i\in G_{g}(\lambda-\eps)\setminus G_g(\lambda+\eps)\Big\}- \pi\big(G_{g}(\lambda-\eps)\setminus G_g(\lambda+\eps)\big)\Big|\rightarrow 0$$
and $\pi\big(G_{g}(\lambda-\eps)\setminus G_g(\lambda+\eps)\big)\leq g_1 \mu_L \big(G_{g}(\lambda-\eps)\setminus G_g(\lambda+\eps)\big),$ 
where $g_1=\max_{x\in S} g(x)$. 

By Lemma \ref{lemdhaus},
\begin{equation*}
\sup_{l_\tau-\frac{\eps_0}{2} \leq \lambda\leq l_\tau+\frac{\eps_0}{2}}\mu_L \big(G_{g}(\lambda-\eps)\setminus G_g(\lambda+\eps)\big)\leq \sup_{l_\tau-\frac{\eps_0}{2} \leq \lambda\leq l_\tau+\frac{\eps_0}{2}}\mu_L \Big( B\big(\partial G_g(\lambda), \frac{3\eps M}{m^2}\big)\Big).
\end{equation*}
For a fixed $\eps>0$, $\mu_L \big( B\big(\partial G_g(\lambda), \frac{3\eps M}{m^2}\big)\big)$  is a continuous function of $\lambda$, and so its maximum is attained in some $\lambda_0\in [l_\tau-\eps_0/2,l_\tau+\eps_0/2]$. Since $reach(\partial G_g(\lambda_0))>0$, the outer Minkowski content of $G_g(\lambda_0)$ and $G_g(\lambda_0)^c$ exist, and so by corollary 3 of \cite{amb08},
$$\sup_{l_\tau-\frac{\eps_0}{2} \leq \lambda\leq l_\tau+\frac{\eps_0}{2}}\mu_L \big(G_{g}(\lambda -\eps)\setminus G_g(\lambda +\eps)\big)=\mathcal{O}(\eps),$$
from which it follows that $\sup_{\lambda \in [l_\tau- \eps_0/2,l_\tau+\eps_0/2]}|\tilde{L}(\lambda)-\hat{L}(\lambda)|\rightarrow 0.$
Using Lemma \ref{lempcont} it follows that 
$\sup_{\lambda \in [l_\tau- \eps_0/2,l_\tau+\eps_0/2]}|L(\lambda)-\tilde{L}(\lambda)|\rightarrow 0,$
then $\sup_{\lambda \in [l_\tau- \eps_0/2, l_\tau+ \eps_0/2]}|L(\lambda)-\hat{L}(\lambda)|\rightarrow 0.$

To prove that $\hat{l}_\tau\rightarrow l_\tau$ a.s., let $0<\eps<\eps_0/2$ and 
$$\gamma=\min\big\{1-\tau-L(l_\tau+\eps/2),L(l_\tau-\eps/2)-(1-\tau)\big\}.$$
 Now observe that $\gamma>0$ since $L$ is decreasing in $l_\tau- \eps_0\leq \lambda\leq l_\tau+\eps_0$. Let $n$ be so large that $\sup_{l_\tau- \eps_0/2\leq \lambda\leq l_\tau+\eps_0/2}|L(\lambda)-\hat{L}(\lambda)|<\gamma/2$. Then $l_\tau-\eps/2<\hat{l}_\tau<l_\tau+\eps/2$.
To conclude the proof, observe that since $\|\nabla g(x)\|>m$ for all $x\in U$,  where $U$ is an open set containing $\overline{G_{g}(l_{\tau}-\eps_0)}\setminus G_g(l_\tau+\eps_0)$, it follows that $\overline{\{x:g(x)<\lambda\}}=\{x:g(x)\leq \lambda\}$ for all $l_\tau-\eps_0<\lambda<l_\tau+\eps_0$. Now we apply theorem 2.1 of \cite{mol98}, which implies that, with probability one,
\begin{equation} \label{eqthunif}
\sup_{l_\tau-\frac{\eps_0}{2}\leq \lambda \leq l_\tau+\frac{\eps_0}{2}} d_H\big(G_{\hat{g}_n}(\lambda),G_g(\lambda)\big)\rightarrow 0.
\end{equation}
Finally the result follows since
$$d_H\big(G_{\hat{g}_n}(\hat{l}_\tau),G_g(l_\tau)\big)\leq d_H\big(G_{\hat{g}_n}(\hat{l}_\tau),G_g(\hat{l}_\tau)\big)+d_H\big(G_g(\hat{l}_\tau),G_g(l_\tau)\big),$$
while \eqref{eqthunif} implies that the first term converges to zero, and the second one converges to zero by Lemma \ref{lemdhaus}.

 \section{Appendix B}

Here we include the proofs of the propositions stated in Section \ref{setup}.

 \textit{Proof of Proposition \ref{prop0}.}\\
 	The proof is based on the ideas used to prove Proposition 1.4 (ii) in \cite{burdzy:06} and the following result (whose proof can be found in \cite{ca:92} 610--613): 
 	$$\inf_{(x,y)\in \overline{D}\times \overline{D}}p(0,x,t,y)=c_t>0,$$
 	where $p(0,x,t,y)$ is the density function introduced in Remark \ref{rem}. Let $A$ be a Borel set such that $\mu_L(A\cap D)>0$. Then for all $t\geq 1$,
 	\begin{equation*}\mathbb{P}_x(T_{A}\leq t)\geq \mathbb{P}_x(T_{A}\leq 1)\geq \int_{A} p(0,x,1,y)dy\geq c_1\mu_L({A\cap C})=c'>0.
 	\end{equation*}
 	By the Markov property, for every $x\in D$, $\mathbb{P}_x(T_{A}\geq k)\leq (1-c')^k$, for all $k\geq 1$, which implies that 
 	$$\sup_{x\in D}\mathbb{E}_x(T_{A})\leq \sup_{x\in D}\sum_{k=0}^\infty \mathbb{P}_x(T_{A}\geq k)<\infty.$$
 	This proves  \eqref{trapcond2} \

  	\textit{Proof of Proposition \ref{properg}}	
 	\begin{proof} 
 		Let $x_0\in D$ and $\eta>0$ be such that $\mathcal{B}(x_0,3\eta)\subset D$. 
 		Since  $\sup_{x\in D} \mathbb E_x T_{\mathcal{B}(x_0,\eta)}<\infty$, 
 		by the Markov inequality there exists an $n_1$ such that 
 		$\inf_{x\in D}\mathbb{P}_x(T_{\mathcal{B}(x_0,\eta)}\leq n_1)>1/2$. 
 		Let $Z_t=x+B_t+\int_0^t \mu(X_s)ds$ be the $d$-dimensional Brownian motion with drift given by $\mu(x)$. 
 		Observe that, since $|\mu(x)|<L$, by Doob's maximal inequality, we have
 		\begin{equation*}
 		\mathbb{P}_x\left(\sup_{s\in [0,t]}|Z_s|<\eta\right)\geq 1-\frac{\sqrt{dt}+Lt}{\eta}.
 		\end{equation*}
 		Now take $t_0$ small enough so that 	
 		$1-\frac{\sqrt{dt_0}+Lt_0}{\eta}=:p_0>0$.	
 		By the strong Markov property,
 		$$\inf_{x\in D}\mathbb{P}_x\big(T_{\mathcal{B}(x_0,\eta)}\leq n_1 \text{ and } X_t\in \mathcal{B}(x_0,2\eta) \text{ for }t\in [T_{\mathcal{B}(x_0,\eta)},T_{\mathcal{B}(x_0,\eta)}+t_0]\big)>\frac{1}{2}p_0.$$
 		Let 
 		$Y=\inf\{n\in \mathbb{N}\colon X_n\in \mathcal{B}(x_0,2\eta)\}$, then $\inf_{x\in D}\mathbb{P}_x(Y\leq n_1+t_0)>p_0/2$. 
 		Applying the Markov property at times $k\lfloor(n_1+t_0)\rfloor$, 
 		$\sup_{x\in D} \mathbb{P}_x(Y\geq k\lfloor(n_1+t_0)\rfloor)\leq (1-p_0/2)^k$, from which it follows that
 		$$
 		\sup_{x\in D}\mathbb{E}_x(Y)\leq \sup_{x\in D}\sum_{k=0}^\infty k\lfloor(n_1+t_0)\rfloor\mathbb{P}_x(Y\geq k\lfloor(n_1+t_0)\rfloor)<\infty.
 		$$
 		Applying theorem 16.0.2 of \cite{MTa}, we obtain, for every $n>0$, that
 		$$
 		\sup_{x\in D}\|\mathbb{P}_x(X_n\in \cdot)-\pi(\cdot)\|_{TV}\leq c_3e^{-c_4n},
 		$$
 		where $c_3,c_4$ are positive finite constants. Using the semigroup property of $\{X_t\}_{t\geq 0}$ 
 		and the fact that $\pi$ is invariant,
 		\begin{align*}
 		\sup_{x\in D}\|\mathbb{P}_x(X_t\in \cdot)-\pi(\cdot)\|_{TV}=&\sup_{x\in D}\left|\int_{D}\mathbb{P}_{y}(X_{t-n}\in \cdot)d\mathbb{P}_x(X_{n}\in dy)-\int_D\mathbb{P}_{y}(X_{t-n}\in \cdot)\pi(y)\right|\\
 		\leq&\sup_{x\in D}\|\mathbb{P}_x(X_n\in \cdot)-\pi(\cdot)\|_{TV},
 		\end{align*}
 		for all $t$ and $n$, with $t\geq n$. \end{proof}
	
 \section{Appendix C}	
 \textit{Proof of Theorem \ref{condrift}}
 \begin{proof}
 Let  $\gamma_n \geq  2h_n$, $\gamma_n \to 0$, $\Delta \to 0$ and  denote  
	$I_n=\{i: X_{t_i}\in B(x,h_n), \exists s_0: t_i< s_0\leq t_{i+1}, X_{s_0}\notin B(x,\gamma_n)\}$.
According to our model, the estimator can be written as
\begin{multline*} 
\hat{\mu}_{n}(x)=\frac{1}{\Delta N_x} \sum_{i=1}^n (B_{t_{i+1}}-B_{t_i})\mathbb{I}_{\{X_{t_i}\in B(x,h_{n})\}}+ \frac{1}{\Delta N_x}\sum_{i\in I_n}\int_{t_i}^{t_{i+1}} \mu(X_s)ds +\\
\frac{1}{\Delta N_x}\sum_{i\in I_n^C}\int_{t_i}^{t_{i+1}} \mu(X_s)ds+ \frac{1}{\Delta N_x} \sum_{ i\in I_n} \int_{t_i}^{t_{i+1}}\eta(X_s)dL_s=:A_{n,T}+B^1_{n,T}+B^2_{n,T}+C_{n,T}.
\end{multline*}
First will prove that $C_{n,T}\rightarrow 0$ in probability. Observe that, we can bound, using Theorem 4.2 in Saisho (1987) 
$$\Big\|\int_{t_i}^{t_{i+1}}\eta(X_s)dL_s\Big\|\leq L_s[t_i,t_{i+1}]\leq C\sqrt{\Delta},$$
being $C$ a positive constant, then
$C_{n,T}\leq C \frac{\# I_n}{\sqrt{\Delta} N_x} \quad a.s$. 
Let us fix $\epsilon>0$, we will prove that 
\begin{equation}\label{boundc}
\mathbb{P}\Big(\frac{\#I_n}{\sqrt{\Delta}N_x}>\epsilon\Big)\rightarrow 0.
\end{equation}

Let $A_{in}= \{\exists s_i: t_i\leq s_i\leq t_{i+1}, X_{s_i}\notin B(x,\gamma_n)\}$. Then,
\begin{multline} \label{bound1}
\mathbb{P}(A_{in}\cap \{X_{t_i}\in B(x,h_n)\} )\leq \\
 \mathbb{P}\Big(\sup_{s\in [t_i,t_{i+1}]} \|X_s-X_{t_i}\|>\gamma_n-h_n|X_{t_i}\in \partial B(x,h_n)\Big)\mathbb{P}(X_{t_i}\in B(x,h_n))\leq \\
2\frac{(\sqrt{2}+\nu)\sqrt{\Delta}}{\gamma_n} \mathbb{P}(X_{t_i}\in B(x,h_n)).
\end{multline}

Consider the random variable $\kappa=\lfloor \epsilon \sqrt{\Delta} N_x\rfloor$. 
Observe that if $\#I_n/(\sqrt{\Delta}N_x)>\epsilon$ then there exists $\{i_1,\ldots,i_\kappa\}$ where $1\leq i_j <n-1$ for all $j=1,\ldots,\kappa$, such that  $\exists s_{i_j}: t_{i_j}< s_{i_j}\leq t_{i_j+1}$ and  $X_{s_{i_j}}\notin B(x,\gamma_n),X_{t_{i_j}}\in B(x,h_n)$ for all $j=1,\ldots,\kappa$. Let us denote $m_n=2(n\epsilon \pi h_n^2f(x)\sqrt{\Delta})$, observe that $m_n\rightarrow \infty$, and from \eqref{bound1} we get
\begin{multline}\label{eq0}
\mathbb{P}\Big(\frac{\#I_n}{\sqrt{\Delta}N_x}>\epsilon\Big)\leq \mathbb{P}\Big(\frac{\#I_n}{\sqrt{\Delta}N_x}>\epsilon, \mathbb{I}_{\{\kappa\leq m_n\}}\Big)+
\mathbb{P}\Big(\kappa> m_n\Big)\\
\leq \sum_{j=1}^{m_n} 2\frac{(\sqrt{2}+\nu)\sqrt{\Delta}}{\gamma_n} \mathbb{P}(X_{t_{i_j}}\in B(x,h_n))+\mathbb{P}\Big(\kappa> m_n\Big).
\end{multline}

By the Ergodic theorem $\kappa/(\epsilon n\pi h_n^2g(x)\sqrt{\Delta})\rightarrow 1$ a.s., then with probability one, for $n$ large enough,  $\kappa\leq m_n$ from where it follows that $\mathbb{P}(\kappa> m_n)\rightarrow 0$. Lastly, again by ergodicity, we have that 
\begin{equation} \label{eqer}
\frac1{m_n\pi h^2_n}\sum_{j=1}^{m_n}\mathbb{P}\left(X_{t_{i_j}}\in B(x,h_{n})\right)\to g(x),
\end{equation}
and  if we choose $\gamma_n$ fulfilling $h_n^4n\Delta/\gamma_n\rightarrow 0$ we get \eqref{boundc} from \eqref{eq0} and \eqref{eqer}.

The proof will be complete if  under our asymptotic scheme, we have
\begin{align}
A_{n,T}&\to 0,\quad \text{ in probability,}\label{a}\\
B^1_{n,T}&\to 0\quad \text{ in probability,}\\
B^2_{n,T}&\to\mu(x)\quad \text{in probability}. \label{b}
\end{align}

Since $\mu$ is Lipschitz 
and $\gamma_n \to 0$, \eqref{b} follows. 
 
Regarding $B^1_{n,T}$ observe that  $\int_{t_i}^{t_{i+1}}\mu(X_s)ds\leq \max_{x \in S}\Vert \mu(x)\Vert \Delta$ and then from \eqref{boundc} we get $B^1_{n,T}\rightarrow 0$ in probability.

Let us consider now \eqref{a}. Each random variable $\mathbb{I}_{\{X_{t_i}\in B(x,h_{n,T})\}}$ is $\mathcal{F}_{t_i}$ measurable,
due to the independence of $B_{t_{i+1}}-B_{t_i}$ w.r.t. $\mathcal{F}_{t_i}$.
Then $E(B_{t_{i+1}}-B_{t_i}|\mathcal{F}_{t_i})=E(B_{t_{i+1}}-B_{t_i})=0$, giving $E(A_{n,T})=0$.
(In fact this proves that the numerator in $A_{n,T}$ is a martingale.)
We now turn to the computation of the variance. First, by the ergodic theorem,
we obtain that
\begin{equation}\label{ergodic2}
{N_x\over n\pi h_n^{2}}\to g(x),\ a.s.
\end{equation}
Defining
$$
\hat{A}_{n,T}=\frac{1}{a_n(x)} \sum_{i=1}^{n-1} (B_{t_{i+1}}-B_{t_i})\mathbb{I}_{\{X_{t_i}\in B(x,h_{n})\}},
$$
with $a_n(x)=\Delta n\pi h^2_ng(x)$, by \eqref{ergodic2}
we know that $A_{n,T}$ and $\hat{A}_{n,T}$ have the same limit in probability. 
Furthermore
\begin{align*}
\mathbb{E}((\hat{A}_{n,T})^2)&=
\frac{1}{a_n(x)^2}\mathbb{E}\left(\sum_{i=1}^{n-1} \mathbb{I}_{\{X_{t_i}\in B(x,h_{n})\}}(B_{t_{i+1}}-B_{t_i})\right)^2\\
&=
\frac{1}{a_n(x)^2} \sum_{i=1}^{n-1}\mathbb{E}\left(\mathbb{I}_{\{X_{t_i}\in B(x,h_{n})\}}(B_{t_{i+1}}-B_{t_i})^2\right) 
\end{align*}
since the cross--terms are zero. 

We then conclude that
\begin{align*}
\mathbb{E}((\hat{A}_{n,T})^2)&=\frac{1}{\left(\Delta n\pi h^2_ng(x)\right)^2}\sum_{i=1}^{n-1}P\left(
\mathbb{I}_{\{X_{t_i}\in B(x,h_{n})\}}
\right)\Delta\\
&\leq 
\frac{1}{\Delta n\pi h^2_ng(x)^2}\frac1{n\pi h^2_n}\sum_{i=1}^{n-1}
\mathbb{P}(X_{t_i}\in B(x,h_{n})).
\end{align*}
By ergodicity, we have
$$
\frac1{n\pi h^2_n}\sum_{i=1}^{n-1}\mathbb{P}\left(X_{t_i}\in B(x,h_{n})\right)\to g(x),
$$
then, taking into account \eqref{ergodic2}, we obtain
$$
\mathbb{E}(({A}_{n,T})^2) \lessapprox \frac{1}{\Delta n\pi h^2_ng(x)} \to 0.$$\end{proof}
  
\section*{Acknowledgements} 
We thank Professor K. Burdzy for helpful comments during the preparation of this manuscript, and two referee's for their constructive comments which improves significantly the present version of the manuscript.
We also thank Dr. Stephen Blake, of the Max Planck Institute for Ornithology, for facilitating access to the data set that was used in this manuscript.
All possible errors are the responsibility of the authors.

\end{document}